\documentclass[a4, 12pt]{amsart}
\usepackage[mathscr]{eucal}
\usepackage{amssymb}
\usepackage{latexsym}
\usepackage{amsthm}
\theoremstyle{plain}

\makeatletter

\title[A lower bound for eigenvalues ]{ A lower bound for
eigenvalues of \\ a clamped plate problem* }
\author{Qing-Ming Cheng and\ Guoxin Wei}
\address{Qing-Ming Cheng\\  Department of Mathematics, Faculty of Science and Engineering,
Saga University, Saga 840-8502,  Japan, cheng@ms.saga-u.ac.jp}
\address{Guoxin Wei\\  Department of Mathematics, Faculty of Science and Engineering,
Saga University, Saga 840-8502,  Japan, wei@ms.saga-u.ac.jp}

\begin{document}
\maketitle

\begin{abstract}
\noindent In this paper, we study eigenvalues of a clamped plate
problem. We obtain a  lower bound for  eigenvalues,  which gives  an
important improvement of  results due to  Levine and Protter.

\end{abstract}



\footnotetext{* Research partially Supported by a Grant-in-Aid for
Scientific Research from JSPS.}

\section *{1. Introduction}

Let $M$ be an $n$-dimensional complete Riemannian manifold.
The following is called {\it Dirichlet eigenvalue problem of Laplacian}:
\begin{equation*}
  {\begin{cases}
     \Delta u =-\lambda u,& \ \ {\rm in} \ \ \Omega ,\\
     u=0 , & \ \ {\rm on}  \ \ \partial \Omega,
     \end{cases}}
     \eqno{(1.1)}
\end{equation*}
where $\Omega$ is a bounded domain in $M$ with piecewise smooth
boundary $\partial\Omega$ and $\Delta$ denotes the Laplacian on $M$.
It is well known that the spectrum of this eigenvalue problem (1.1)
is real and discrete.
$$0<\lambda_1<\lambda_2\leq\lambda_3\leq\cdots\rightarrow\infty,$$
where each $\lambda_i$ has finite multiplicity which is repeated
according to its multiplicity.

Let $V(\Omega)$ denotes the volume of $\Omega$ and $B_n$ the volume
of the unit ball in $\mathbf{R}^n$, then the following Weyl's
asymptotic formula holds
$$\lambda_k\sim\frac{4\pi^2}{(B_nV(\Omega))^{\frac{2}{n}}}k^{\frac{2}{n}},
\ \ \ k\rightarrow\infty.\eqno{(1.2)}$$ From this asymptotic
formula, one can infer
$$\frac{1}{k}\sum_{i=1}^k\lambda_i\sim\frac{n}{n+2}\frac{4\pi^2}{(B_nV(\Omega))^{\frac{2}{n}}}k^{\frac{2}{n}},
\ \ \ k\rightarrow\infty.\eqno{(1.3)}$$ In particular, when
$M=\mathbf{R}^n$, P\'olya \cite{[P]} proved
$$\lambda_k\geq\frac{4\pi^2}{(B_nV(\Omega))^{\frac{2}{n}}}k^{\frac{2}{n}},\
\ \ \rm{for}\  k=1,2,\cdots,\eqno{(1.4)}$$ if $\Omega$ is a tiling
domain in $\mathbf{R}^n$. Moreover,  he conjectured for a general
bounded domain,

\par\bigskip \noindent{\bf Conjecture of P\'olya}. {\it
If $\Omega$ is a bounded domain in $\mathbf{R}^n$, then eigenvalue
$\lambda_k$ of the eigenvalue problem $(1.1)$ satisfies
$$\lambda_k\geq\frac{4\pi^2}{(B_nV(\Omega))^{\frac{2}{n}}}k^{\frac{2}{n}},\
\ \ \rm{for}\ k=1,2,\cdots.\eqno{(1.5)}$$} On the conjecture of
P\'olya, Li and Yau \cite{[LY]} (cf. \cite{[B]}, \cite{[L]}) proved
$$\frac{1}{k}\sum_{i=1}^k\lambda_i\geq\frac{n}{n+2}\frac{4\pi^2}{(B_nV(\Omega))^{\frac{2}{n}}}k^{\frac{2}{n}},\
\ \ \rm{for}\ k=1,2,\cdots.\eqno{(1.6)}$$
It is sharp about the
highest order term of $k$ in the sense of average according to
(1.3). From this formula, one can derive
$$\lambda_k\geq\frac{n}{n+2}\frac{4\pi^2}{(B_nV(\Omega))^{\frac{2}{n}}}k^{\frac{2}{n}},\
\ \ \rm{for}\ k=1,2,\cdots,\eqno{(1.7)}$$ which gives a partial
solution for the conjecture of P\'olya with a factor
$\frac{n}{n+2}$.

Furthermore, Melas \cite{[M]} obtained the following estimate which is
an improvement of (1.6).
$$\frac{1}{k}\sum_{i=1}^k\lambda_i\geq\frac{n}{n+2}\frac{4\pi^2}{(B_nV(\Omega))^{\frac{2}{n}}}k^{\frac{2}{n}}
+c_n\frac{V(\Omega)}{I(\Omega)},\ \ \ \rm{for}\
k=1,2,\cdots,\eqno{(1.8)}$$ where  $c_n$ is a constant depending only on
the dimension $n$ and
$$
I(\Omega)=\min\limits_{a\in \mathbf{R}^n}\int_\Omega|x-a|^2dx
$$
 is called {\it  the moment of inertia} of
$\Omega$.

For a bounded domain in an $n$-dimensional complete Riemannian manifold,  Cheng and Yang \cite{[CY2]}
have also given a lower bound for eigenvalues, recently.

Our purpose in this paper is to study eigenvalues of  the following  clamped plate problem.
Let $\Omega$ be a bounded
domain in an $n$-dimensional complete Riemannian manifold $M^n$.
The following is called {\it a clamped plate problem}, which describes characteristic vibrations
of a clamped plate:

\begin{equation*}
  {\begin{cases}
     \Delta^2 u = \Gamma  u,& \ \ {\rm in} \ \ \Omega ,\\
     u=\dfrac{\partial u}{\partial \nu}=0 , & \ \ {\rm on}  \ \ \partial \Omega,
     \end{cases}}
     \eqno{(1.9)}
\end{equation*}
where $\Delta$ is the Laplacian   on $M^n$ and
$\nu$ denotes the outward unit normal to the boundary $\partial \Omega$. It is well known that this
problem has a real and discrete spectrum

$$0<\Gamma_1\leq\Gamma_2\leq\cdots\leq\Gamma_k\leq\cdots\to
+\infty,$$ where each $\Gamma_i$ has finite multiplicity which is
repeated according to its multiplicity.

For eigenvalues of the clamped plate problem,  Agmon  \cite{[Ag]} and Pleijel \cite{[P1]}
gave  the following   asymptotic  formula,  which is a generalization of Weyl's asymptotic formula,
$$
\Gamma_k\sim
\dfrac{16\pi^4}{\big(B_nV(\Omega)\big)^{\frac{4}{n}}}k^{\frac{4}{n}},\
\ \  k\rightarrow\infty.
  \eqno{(1.10)}
  $$
The average of the eigenvalues satisfies
$$
\frac{1}{k}\sum_{j=1}^k\Gamma_j
\sim\frac{n}{n+4}\dfrac{16\pi^4}{\big(B_nV(\Omega)\big)^{\frac{4}{n}}}k^{\frac{4}{n}},
\ \ k\rightarrow\infty.\eqno{(1.11)}
$$
Furthermore, Levine and Protter \cite{[LP]} proved that eigenvalues of
the clamped plate problem satisfy
$$\dfrac{1}{k}\sum_{j=1}^k\Gamma_j
\geq\frac{n}{n+4}\dfrac{16\pi^4}{\big(B_nV(\Omega)\big)^{\frac{4}{n}}}k^{\frac{4}{n}}.\eqno{(1.12)}$$
The inequality  (1.12) is sharp about the highest order term of $k$
according to (1.11).

 In this paper, we give an important improvement of the result due to Levine and Protter \cite{[LP]}
 by adding to its right hand side two terms of the lower order terms of $k$.  In fact, we prove the
following:

\par\bigskip \noindent{\bf Theorem}. {\it Let $\Omega$ be a
bounded domain in an $n$-dimensional Euclidean space $\mathbf{R}^n$.
Eigenvalues of the clamped plate problem satisfy
$$
\aligned \frac{1}{k}\sum_{j=1}^k\Gamma_j\geq&
\frac{n}{n+4}\dfrac{16\pi^4}{\big(B_nV(\Omega)\big)^{\frac{4}{n}}}k^{\frac{4}{n}}\\
 &+\left(\frac{n+2}{12n(n+4)}-\frac{1}{1152n^2(n+4)}
    \right)\frac{V(\Omega)}{I(\Omega)}\frac{n}{n+2}
 \dfrac{4\pi^2}{\big(B_nV(\Omega)\big)^{\frac{2}{n}}}k^{\frac{2}{n}}\\
 &+\left(\frac{1}{576n(n+4)}-\frac{1}{27648n^2(n+2)(n+4)}\right)
 \left(\frac{V(\Omega)}{I(\Omega)}\right)^2,
\endaligned
\eqno{(1.13)}
$$
where $I(\Omega)$ is the moment of inertia  of $\Omega$. }

\par\bigskip \noindent{\bf Corollary}. {\it Let $\Omega$ be a
bounded domain in an $n$-dimensional Euclidean space $\mathbf{R}^n$.
Then eigenvalues $\Gamma_j$'s of the clamped plate problem satisfy
$$
\aligned \frac{1}{k}\sum_{j=1}^k\Gamma_j\geq&
\frac{n}{n+4}\dfrac{16\pi^4}{\big(B_nV(\Omega)\big)^{\frac{4}{n}}}k^{\frac{4}{n}}\\
 &+\left(\frac{n+2}{12n(n+4)}-\frac{1}{1152n^2(n+4)}
    \right)\frac{1}{\sum_{i=1}^n\mu_i^{-1}}\frac{n}{n+2}
 \dfrac{4\pi^2}{\big(B_nV(\Omega)\big)^{\frac{2}{n}}}k^{\frac{2}{n}}\\
 &+\left(\frac{1}{576n(n+4)}-\frac{1}{27648n^2(n+2)(n+4)}\right)
 \left(\frac{1}{\sum_{i=1}^n\mu_i^{-1}}\right)^2,
\endaligned
\eqno{(1.14)}
$$
where $\mu_1$, $\cdots$, $\mu_n$ are the first $n$ nonzero
eigenvalues of the Neumann eigenvalue problem of Laplacian
\begin{equation*}
  {\begin{cases}
     \Delta v = -\mu  v,& \ \ {\rm in} \ \ \Omega ,\\
     \dfrac{\partial v}{\partial \nu}=0 , & \ \ {\rm on}  \ \ \partial \Omega.
     \end{cases}}
\end{equation*}}


\par\bigskip \noindent{\bf Remark 1}. On universal estimates for
eigenvalues of the clamped plate problem, one can see \cite{[CHW]},
\cite{[CIM2]}, \cite{[CY1]}, \cite{[H]}  and  \cite{[WX2]}.

\section*{2. Proof of results}
For a bounded domain $\Omega$, {\it the moment of inertia} of
$\Omega$ is defined by
$$
I(\Omega)=\min\limits_{a\in \mathbf{R}^n}\int_{\Omega}|x-a|^2dx.
$$
By a translation of the origin and a suitable rotation of axes,
we can assume that the
center of mass is the origin and
$$
I(\Omega)=\int_\Omega|x|^2dx.
$$

For reader's convenience, we first review the definition and serval
properties of  the symmetric decreasing rearrangements. Let
$\Omega\subset \mathbf{R}^n$ be a bounded domain. Its {\it symmetric
rearrangement} $\Omega^*$ is the open ball with the same volume as
$\Omega$,
$$\Omega^*=\big\{x\in
\mathbf{R}^n| \ \   |x|<\biggl(\dfrac{{\rm Vol}
(\Omega)}{B_n}\biggl)^{\frac1n}\big\}.$$ By using a symmetric
rearrangement of $\Omega$, we have
$$
I(\Omega)=\int_\Omega |x|^2dx\geq\int_{\Omega^*}
|x|^2dx=\frac{n}{n+2}V(\Omega)\left(\frac{V(\Omega)}{B_n}\right)^{\frac{2}{n}}.\eqno{(2.1)}
$$

 Let $h$ be a nonnegative bounded continuous
function on $\Omega$, we can consider its {\it distribution
function} $\mu_h(t)$ defined by
$$
\mu_h(t)={\rm Vol}(\{x\in\Omega| \ \ h(x)>t\}).
$$
The distribution function can be viewed as a function from
$[0,\infty)$ to $[0, V(\Omega)]$. The {\it symmetric decreasing
rearrangement $h^*$ of $h$} is defined by
$$h^*(x)=\inf \{t\geq0|\mu_h(t)<B_n|x|^n\}$$
for $x\in \Omega^*$. By definition, we know that ${\rm Vol}(\{x\in
\Omega|h(x)>t\})={\rm Vol}(\{x\in \Omega^*|h^*(x)>t\}), \ \forall
t>0$ and $h^*(x)$ is a radially symmetric function.

Putting $g(|x|):=h^*(x)$, one gets that $g: [0, +\infty)\rightarrow
[0, \sup {h}]$ is a non-increasing function of $|x|$.
 Using the well known properties of the symmetric decreasing
rearrangement, we obtain
$$\int_{\mathbf{R}^n}h(x)dx=\int_{\mathbf{R}^n}h^*(x)dx=nB_n\int_0^\infty
s^{n-1}g(s)ds\eqno{(2.2)}
$$
and
$$\int_{\mathbf{R}^n}|x|^4h(x)dx\geq\int_{\mathbf{R}^n}|x|^4h^*(x)dx=nB_n\int_0^\infty
s^{n+3}g(s)ds.\eqno{(2.3)}
$$
Good sources of further information on rearrangements are
\cite{[Ba]}, \cite{[Po]}.

One gets from the coarea formula that
$$
\mu_h(t)=\int_t^{\sup h}\int_{\{h=s\}}|\nabla h|^{-1}d\sigma_sds.
$$
Since $h^*$ is radial, we have  
\begin{equation*}
\begin{aligned}
\mu_h(g(s))&={\rm Vol}\{x\in
\Omega|h(x)>g(s)\}={\rm Vol}\{x\in \Omega^*|h^*(x)>g(s)\}\\
&={\rm Vol}\{x\in \Omega^*|g(|x|)>g(s)\}=B_ns^n.
\end{aligned}
\end{equation*}
 It follows that
$$
nB_ns^{n-1}=\mu_h^{'}(g(s))g^{'}(s)
$$
for almost every $s$. Putting $\tau:=\sup |\nabla h|$, we obtain
from the above equations and the isoperimetric inequality that
$$
\aligned -\mu_h^{'}(g(s))&=\int_{\{h=g(s)\}}|\nabla
   h|^{-1}d\sigma_{g(s)}\geq\tau^{-1}{\rm Vol}_{n-1}(\{h=g(s)\})\\
   &\geq\tau^{-1}nB_ns^{n-1}.
\endaligned
$$
Therefore, one  obtains
$$-\tau\leq g^{'}(s)\leq 0\eqno{(2.4)}$$
for almost every $s$.

Next, we prepare  the following lemma in order to prove of our
theorem.

\par\bigskip \noindent{\bf Lemma 2.1}. {\it Let $b\geq1$, $\eta>0$
and $\psi: [0, +\infty)\rightarrow [0, +\infty)$ be a
decreasing smooth function such that
$$-\eta\leq \psi^{'}(s)\leq 0$$
and, for a constant $d<1$,
$$
\dfrac{\psi(0)^{\frac{2b+2}{b}}}{6b\eta^2(bA)^{\frac{2}{b}}}<d
$$
with
$$
A:=\int_0^{\infty}s^{b-1}\psi(s)ds>0.
$$
Then, we have
$$
\aligned \int_0^\infty
s^{b+3}\psi(s)ds\geq&\dfrac{1}{b+4}(bA)^{\frac{b+4}{b}}\psi(0)^{-\frac{4}{b}}\\
    &+\left(\dfrac{1}{3b(b+4)\eta^2}-\dfrac{d}{6(b+2)^2(b+4)\eta^2}
    \right)(bA)^{\frac{b+2}{b}}\psi(0)^{\frac{2b-2}{b}}\\
   &+\left(\dfrac{1}{36b(b+4)\eta^4}-\dfrac{d}{36(b+2)^2(b+4)\eta^4}\right)A\psi(0)^4.
\endaligned
\eqno{(2.5)}
$$
}

\vskip 3pt \noindent {\it Proof.}  Defining
$$
D:=\int_0^\infty s^{b+1}\psi(s)ds,
$$
one can prove from the same assertions as in the lemma 1 of \cite{[M]},
$$
D=\int_0^\infty s^{b+1}\psi(s)ds\geq
\frac{1}{b+2}(bA)^{\frac{b+2}{b}}\psi(0)^{-\frac{2}{b}}+\frac{A\psi(0)^2}{6(b+2)\eta^2}.\eqno{(2.6)}
$$
Since the formula (2.6) holds for any constant $b\geq 1$, we have
\begin{equation*}
\begin{aligned}
 &\int_0^{\infty} s^{b+3}\psi(s)ds\\
&\geq\frac{1}{b+4}((b+2)D)^{\frac{b+4}{b+2}}\psi(0)^{-\frac{2}{b+2}}+\frac{D\psi(0)^2}{6(b+4)\eta^2}\\
    &\geq\frac{1}{b+4}\biggl[(bA)^{\frac{b+2}{b}}\psi(0)^{-\frac{2}{b}}
    +\frac{A\psi(0)^2}{6\eta^2}\biggl]^{\frac{b+4}{b+2}}
    \psi(0)^{-\frac{2}{b+2}}\\
    &\ \ \ +\frac{\psi(0)^2}{6(b+4)\eta^2}\biggl[\frac{1}{b+2}(bA)^{\frac{b+2}{b}}\psi(0)^{-\frac{2}{b}}
    +\frac{A\psi(0)^2}{6(b+2)\eta^2}\biggl]\\
     &=\frac{1}{b+4}\biggl[(bA)^{\frac{b+2}{b}}\psi(0)^{-\frac{2}{b}}+\frac{A\psi(0)^2}{6\eta^2}\biggl]
    \biggl[(bA)^{\frac{b+2}{b}}\psi(0)^{-\frac{2}{b}}\biggl]^{\frac{2}{b+2}}\\
    &\ \ \times
   \biggl(1+\frac{A\psi(0)^{\frac{2b
   +2}{b}}}{6(bA)^{\frac{b+2}{b}}\eta^2}\biggl)^{\frac{2}{b+2}}\psi(0)^{-\frac{2}{b+2}}\\
    &\ \ \ +\frac{1}{6(b+2)(b+4)\eta^2}(bA)^{\frac{b+2}{b}}\psi(0)^{\frac{2b-2}{b}}
     +\frac{A\psi(0)^4}{36(b+2)(b+4)\eta^4}\\
     &\geq\frac{1}{b+4}\biggl[(bA)^{\frac{b+2}{b}}\psi(0)^{-\frac{2}{b}}+\frac{A\psi(0)^2}{6\eta^2}\biggl]
    \biggl[(bA)^{\frac{b+2}{b}}\psi(0)^{-\frac{2}{b}}\biggl]^{\frac{2}{b+2}}\\
    &\ \ \ \times\left\{1+\frac{1}{b+2}\frac{A\psi(0)^{\frac{2b+2}{b}}}{6(bA)^{\frac{b+2}{b}}\eta^2}
    \biggl(2-\frac{b}{b+2}\frac{A\psi(0)^{\frac{2b+2}{b}}}{6(bA)^{\frac{b+2}{b}}\eta^2}\biggl)\right\}
    \psi(0)^{-\frac{2}{b+2}}\\
    &\qquad\qquad\qquad (\text{from the Taylor formula})\\
    &\ \ \ +\frac{1}{6(b+2)(b+4)\eta^2}(bA)^{\frac{b+2}{b}}\psi(0)^{\frac{2b-2}{b}}
     +\frac{A\psi(0)^4}{36(b+2)(b+4)\eta^4}\\
 &\geq\frac{1}{b+4}\biggl[(bA)^{\frac{b+2}{b}}\psi(0)^{-\frac{2}{b}}+\frac{A\psi(0)^2}{6\eta^2}\biggl]
    \biggl[(bA)^{\frac{b+2}{b}}\psi(0)^{-\frac{2}{b}}\biggl]^{\frac{2}{b+2}}\\
    &\ \ \ \times\left\{1+\frac{1}{b+2}\frac{A\psi(0)^{\frac{2b+2}{b}}}{6(bA)^{\frac{b+2}{b}}\eta^2}
   \biggl(2-\frac{b}{b+2}d\biggl)\right\}
    \psi(0)^{-\frac{2}{b+2}}\\
    &\ \ \ +\frac{1}{6(b+2)(b+4)\eta^2}(bA)^{\frac{b+2}{b}}\psi(0)^{\frac{2b-2}{b}}
     +\frac{A\psi(0)^4}{36(b+2)(b+4)\eta^4}\\
        &= \frac{1}{b+4}(bA)^{\frac{b+4}{b}}\psi(0)^{-\frac{4}{b}}\\
    &\ \ \ +\left(\frac{1}{3b(b+4)\eta^2}-\frac{d}{6(b+2)^2(b+4)\eta^2}
    \right)(bA)^{\frac{b+2}{b}}\psi(0)^{\frac{2b-2}{b}}\\
   &\ \ \
   +\left(\frac{1}{36b(b+4)\eta^4}-\frac{d}{36(b+2)^2(b+4)\eta^4}\right)A\psi(0)^4.
\end{aligned}
\end{equation*}
This completes the proof of the lemma.
$$\eqno{\Box}$$

\vskip 3pt\noindent {\it Proof of  Theorem}.  Let $u_j$ be an orthonormal eigenfunction
corresponding to the eigenvalue $\Gamma_j$, that is, $u_j$ satisfies
\begin{equation*}
  {\begin{cases}
     \Delta^2 u_j = \Gamma_j  u_j,& \ \ {\rm in} \ \ \Omega ,\\
     u_j=\dfrac{\partial u_j}{\partial \nu}=0 , & \ \ {\rm on}  \ \ \partial \Omega, \\
     \int u_iu_j=\delta_{ij}, & \ \ \text{for any $i$, $j$}.
     \end{cases}}
     \eqno{(2.7)}
\end{equation*}
Thus, $\{u_j\}_{j=1}^{\infty}$ forms an orthonormal basis of $L^2(\Omega)$.
We define a function $\varphi_j$ by
\begin{equation*}
   \varphi_j(x)={\begin{cases}
    u_j(x), & \ \  x\in \Omega ,\\
     0 , & \ \ x\in \mathbf{R}^n\setminus\Omega. \\
         \end{cases}}
\end{equation*}
Denote by $\widehat{\varphi}_j(z)$ the
Fourier transform of $\varphi_j(x)$. For any $z\in \mathbf{R}^n$, we have
by definition that
$$\widehat{\varphi}_j(z)=(2\pi)^{-n/2}\int_{\mathbf{R}^n} \varphi_j(x)e^{i<x,z>}dx
=(2\pi)^{-n/2}\int_\Omega u_j(x)e^{i<x,z>}dx.\eqno{(2.8)}$$
 From the Plancherel formula, we have
 $$
 \int_{\mathbf{R}^n}{\widehat{\varphi}_i}(z){\widehat{\varphi}_j}(z)dz=\delta_{ij}
 $$
for any $i, j$.
Since $\{u_j\}_{j=1}^{\infty}$ is an orthonormal basis in $L^2(\Omega)$,
the Bessel inequality implies that
$$\sum_{j=1}^k|\widehat{\varphi}_j(z)|^2\leq
(2\pi)^{-n}\int_\Omega|e^{i<x,z>}|^2dx=(2\pi)^{-n}V(\Omega).\eqno{(2.9)}$$
For each $q=1, \cdots, n$, $j=1, \cdots, k$, we deduce from the
divergence theorem and  $u_j|_{\partial\Omega}=\dfrac{\partial
u_j}{\partial \nu}|_{\partial\Omega}=0$ that
$$
\aligned
z_q^2\widehat{\varphi}_j(z)&=(2\pi)^{-n/2}\int_{\mathbf{R}^n}\varphi_j(x)(-i)^2\frac{\partial^2
e^{i<x,z>}}{\partial x_q^2}dx\\
&=-(2\pi)^{-n/2}\int_{\mathbf{R}^n}\frac{\partial^2
\varphi_j(x)}{\partial
x_q^2}e^{i<x,z>}dx\\
&=-\widehat{\frac{\partial^2 \varphi_j}{\partial x_q^2}}(z).
\endaligned
\eqno{(2.10)}
$$
It follows from the Parseval's identity   that
$$
\aligned
\int_{\mathbf{R}^n}|z|^4|\widehat{\varphi}_j(z)|^2dz&=\int_{\mathbf{R}^n}\left
       ||z|^2\widehat{\varphi}_j(z)\right|^2dz\\
      &=\int_{\mathbf{R}^n}|\sum_{q=1}^n\widehat{\frac{\partial^2
      \varphi_j}{\partial x_q^2}}(z)|^2dz\\
      &=\int_\Omega(\sum_{q=1}^n\frac{\partial^2 u_j}{\partial
      x_q^2})^2dx\\
      &=\int_\Omega|\Delta u_j(x)|^2dx\\
      &=\int_\Omega  u_j(x)\Delta^2 u_j(x)dx\\
      &=\int_\Omega \Gamma_ju_j^2(x)dx\\
      &=\Gamma_j.
\endaligned
\eqno{(2.11)}
$$
Since
$$\nabla \widehat{\varphi}_j(z)=(2\pi)^{-n/2}\int_\Omega ixu_j(x)e^{i<x,z>}dx,\eqno{(2.12)}$$
we obtain
$$\sum_{j=1}^k|\nabla \widehat{\varphi}_j(z)|^2\leq
(2\pi)^{-n}\int_\Omega|ixe^{i<x,z>}|^2dx=(2\pi)^{-n}I(\Omega).\eqno{(2.13)}$$
Putting
$$
h(z):=\sum_{j=1}^k|\widehat{\varphi}_j(z)|^2,
$$
one derives from (2.9) that $0\leq h(z)\leq (2\pi)^{-n}V(\Omega),$
it follows from (2.13) and the Cauchy-Schwarz inequality that
$$
\aligned |\nabla h(z)|&\leq
2\left(\sum_{j=1}^k|\widehat{\varphi}_j(z)|^2\right)^{1/2}\left(\sum_{j=1}^k|\nabla
     \widehat{\varphi}_j(z)|^2\right)^{1/2}\\
     &\leq2(2\pi)^{-n}\sqrt{V(\Omega)I(\Omega)}
\endaligned
\eqno{(2.14)}
$$
for every $z\in \mathbf{R}^n$. From the Parseval's identity, we
derive
$$
\int_{\mathbf{R}^n}h(z)dz=\sum_{j=1}^k\int_\Omega
|u_j(x)|^2dx=k.\eqno{(2.15)}
$$
Applying the symmetric decreasing rearrangement to $h$ and noting
that $\tau=\sup |\nabla
h|\leq2(2\pi)^{-n}\sqrt{V(\Omega)I(\Omega)}:=\eta$, we obtain, from
(2.4),
$$
-\eta\leq-\tau\leq g^{'}(s)\leq 0\eqno{(2.16)}
$$
for almost every $s$. According to (2.2) and (2.15), we infer
$$k=\int_{\mathbf{R}^n}h(z)dz=\int_{\mathbf{R}^n}h^*(z)dz=nB_n\int_0^\infty
s^{n-1}g(s)ds.\eqno{(2.17)}
$$
From (2.3) and (2.11), we obtain
$$
\aligned
\sum_{j=1}^k\Gamma_j&=\int_{\mathbf{R}^n}|z|^4h(z)dz\\
     &\geq\int_{\mathbf{R}^n}|z|^4h^*(z)dz\\
     &=nB_n\int_0^\infty s^{n+3}g(s)ds.
\endaligned
\eqno{(2.18)}
$$
In order to apply Lemma 2.1, from (2.17) and the definition of $A$ ,
we take
 $$
 \psi(s)=g(s),\ \
 A=\frac{k}{nB_n},\ \
\eta=2(2\pi)^{-n}\sqrt{V(\Omega)I(\Omega)},\eqno{(2.19)}
$$
from (2.1), we deduce  that
$$
\eta\geq2(2\pi)^{-n}\left(\frac{n}{n+2}\right)^{\frac{1}{2}}B_n^{-\frac{1}{n}}V(\Omega)^{\frac{n+1}{n}}
.\eqno{(2.20)}
$$ 
On the other hand, 
 $0<g(0)\leq\sup h^*(z)=\sup h(z)\leq(2\pi)^{-n}V(\Omega),$
 we have from (2.1), (2.19) and (2.20) that
$$
\aligned &\ \ \
\frac{g(0)^{\frac{2n+2}{n}}}{6n\eta^2(nA)^{\frac{2}{n}}}
          \leq\frac{((2\pi)^{-n}V(\Omega))^{\frac{2n+2}{n}}}
          {6n(2(2\pi)^{-n}\left(\frac{n}{n+2}\right)^{\frac{1}{2}}B_n^{-\frac{1}{n}}V(\Omega)^{\frac{n+1}{n}})^2
          (\frac{k}{B_n})^{\frac{2}{n}}}\\
          &=\frac{n+2}{24n^2}\frac{B_n^{\frac{4}{n}}}{(2\pi)^2k^{\frac{2}{n}}}
           \leq \frac{n+2}{24n^2}\frac{B_n^{\frac{4}{n}}}{(2\pi)^2}.
\endaligned
$$
 By a direct
calculation, one sees from
$B_n=\dfrac{2\pi^{\frac{n}{2}}}{n\Gamma(\frac{n}{2})}$ that
$$
\frac{B_n^{\frac{4}{n}}}{(2\pi)^2}<\frac{1}{2},\eqno{(2.21)}
$$
where $\Gamma(\frac n2)$ is the Gamma function. From the above arguments, one
has
$$
\frac{g(0)^{\frac{2n+2}{n}}}{6n\eta^2(nA)^{\frac{2}{n}}}\leq
\frac{n+2}{48n^2}<1.\eqno{(2.22)}
$$
Hence we know that the function $\psi(s)=g(s)$ satisfies the conditions in Lemma 2.1
with $b=n$ and
$$
\eta=2(2\pi)^{-n}\sqrt{V(\Omega)I(\Omega)},\ \ d=\frac{n+2}{48n^2}.
$$
 From  Lemma 2.1 and (2.18), we conclude
$$
\aligned
\sum_{j=1}^k\Gamma_j&\geq nB_n\int_0^\infty s^{n+3}g(s)ds\\
    &\geq \frac{n}{n+4}(B_n)^{-\frac{4}{n}}k^{\frac{n+4}{n}}g(0)^{-\frac{4}{n}}\\
    &\ \ +\left(\frac{1}{3(n+4)\eta^2}-\frac{1}{288n(n+2)(n+4)\eta^2}
    \right)k^{\frac{n+2}{n}}(B_n)^{-\frac{2}{n}}g(0)^{\frac{2n-2}{n}}\\
   &+\left(\frac{1}{36n(n+4)\eta^4}-\frac{1}{1728n^2(n+2)(n+4)\eta^4}\right)k g(0)^4.
\endaligned
\eqno{(2.23)}
$$
Defining a   function $F$ by
$$
\aligned
F(t)=&\frac{n}{n+4}(B_n)^{-\frac{4}{n}}k^{\frac{n+4}{n}}t^{-\frac{4}{n}}\\
   &\ \ +
   \left(\frac{1}{3(n+4)\eta^2}-\frac{1}{288n(n+2)(n+4)\eta^2}
    \right)k^{\frac{n+2}{n}}(B_n)^{-\frac{2}{n}}t^{\frac{2n-2}{n}}\\
    & \ \ +\left(\frac{1}{36n(n+4)\eta^4}-\frac{1}{1728n^2(n+2)(n+4)\eta^4}\right)k
    t^4.
\endaligned
\eqno{(2.24)}
$$
It is not hard to prove from (2.20) that $\eta
\geq(2\pi)^{-n}B_n^{-\frac{1}{n}}V(\Omega)^{\frac{n+1}{n}}$. 
Furthermore,  it follows from (2.24) that
$$
\aligned
&F^{'}(t)\\
&\leq-\frac{4}{n+4}(B_n)^{-\frac{4}{n}}k^{\frac{n+4}{n}}t^{-1-\frac{4}{n}}\\
   &\ \ \ +\left(\frac{2(n-1)}{3n(n+4)}-\frac{(n-1)}{144n^2(n+2)(n+4)}
    \right)k^{\frac{n+2}{n}}(2\pi)^{2n}V(\Omega)^{-\frac{2(n+1)}{n}} t^{\frac{n-2}{n}}\\
    &\ \ \ +\left(\frac{1}{9n(n+4)}-\frac{1}{432n^2(n+2)(n+4)}\right)k
    t^3(2\pi)^{4n}(B_n)^{\frac{4}{n}}V(\Omega)^{-\frac{4(n+1)}{n}}\\
    &=\frac{k}{n+4}t^{-\frac{n+4}{n}}\times
    \big\{(\frac{2(n-1)}{3n}-\frac{(n-1)}{144n^2(n+2)}
    )(2\pi)^{2n}k^{\frac{2}{n}}
    V(\Omega)^{-\frac{2(n+1)}{n}}t^{\frac{2n+2}{n}}\\
    &\ \ \ -4(B_n)^{-\frac{4}{n}}k^{\frac{4}{n}}
    +(\frac{1}{9n}-\frac{1}{432n^2(n+2)})
      (2\pi)^{4n}(B_n)^{\frac{4}{n}}V(\Omega)^{-\frac{4(n+1)}{n}}t^{\frac{4n+4}{n}}\big\}.
\endaligned
$$
Hence, we have
$$
\aligned
&\frac{n+4}{k}t^{\frac{n+4}{n}}F^{'}(t)\\
&\leq
    \big(\frac{2(n-1)}{3n}-\frac{(n-1)}{144n^2(n+2)}
    \big)(2\pi)^{2n}k^{\frac{2}{n}}
    V(\Omega)^{-\frac{2(n+1)}{n}}t^{\frac{2n+2}{n}}\\
    &\ \ \ -4(B_n)^{-\frac{4}{n}}k^{\frac{4}{n}}
    +\big(\frac{1}{9n}-\frac{1}{432n^2(n+2)}\big)
      (2\pi)^{4n}(B_n)^{\frac{4}{n}}V(\Omega)^{-\frac{4(n+1)}{n}}t^{\frac{4n+4}{n}}.
\endaligned
\eqno{(2.25)}
$$
Since the right hand side of (2.25) is an increasing function of
$t$, if it is not larger  than $0$ at $t=(2\pi)^{-n}V(\Omega)$, that
is,
$$
\aligned &\big(\frac{2(n-1)}{3n}-\frac{(n-1)}{144n^2(n+2)}
    \big)(2\pi)^{2n}k^{\frac{2}{n}}
    V(\Omega)^{-\frac{2(n+1)}{n}}((2\pi)^{-n}V(\Omega))^{\frac{2n+2}{n}}\\
    &\ \ \
    +\big(\frac{1}{9n}-\frac{1}{432n^2(n+2)}\big)
      (2\pi)^{4n}(B_n)^{\frac{4}{n}}V(\Omega)^{-\frac{4(n+1)}{n}}((2\pi)^{-n}V(\Omega))^{\frac{4n+4}{n}}\\
      &\ \ \ -4(B_n)^{-\frac{4}{n}}k^{\frac{4}{n}}\leq0,
\endaligned
\eqno{(2.26)}
$$
then one has from (2.25) that $F^{'}(t)\leq0$ on $(0,
(2\pi)^{-n}V(\Omega)]$. Hence, $F(t)$ is decreasing on $(0,
(2\pi)^{-n}V(\Omega)]$. Indeed, by a direct calculation, we have
that (2.26) is equivalent to
$$
\aligned &\ \
\left(\frac{(n-1)}{6n}-\frac{(n-1)}{576n^2(n+2)}\right)(2\pi)^{-2}k^{\frac{2}{n}}\\
&\ \
+\left(\frac{1}{36n}-\frac{1}{1728n^2(n+2)}\right)(2\pi)^{-4}(B_n)^{\frac{4}{n}}
   \\
    &\leq (B_n)^{-\frac{4}{n}}k^{\frac{4}{n}}.
\endaligned
\eqno{(2.27)}
$$
From (2.21), we can prove that $(2\pi)^{-2}(B_n)^{\frac{4}{n}}<1$
and
$$
\aligned &\ \
   \left(\frac{(n-1)}{6n}-\frac{(n-1)}{576n^2(n+2)}\right)(2\pi)^{-2}k^{\frac{2}{n}}\\
    &+\left(\frac{1}{36n}-\frac{1}{1728n^2(n+2)}\right)(2\pi)^{-4}(B_n)^{\frac{4}{n}}\\
&<\frac{1}{6}(2\pi)^{-2}k^{\frac{2}{n}}+\frac{1}{36n}(2\pi)^{-2}\\
   &<(2\pi)^{-2}\left\{\frac{1}{6}k^{\frac{4}{n}}+\frac{1}{36n}\right\}\\
   &<(2\pi)^{-2}k^{\frac{4}{n}}< (B_n)^{-\frac{4}{n}}k^{\frac{4}{n}},
\endaligned
\eqno{(2.28)}
$$
that is, $F(t)$ is a decreasing function on $(0, (2\pi)^{-n}V(\Omega)]$.

On the other hand, since $0<g(0)\leq(2\pi)^{-n}V(\Omega)$ and the
right hand side of the formula (2.23) is $F(g(0))$, which is a
decreasing function of $g(0)$ on $(0, (2\pi)^{-n}V(\Omega)]$, then
we can replace $g(0)$ by $(2\pi)^{-n}V(\Omega)$ in (2.23) which
gives inequality
$$
\aligned \frac{1}{k}\sum_{j=1}^k\Gamma_j\geq&
\frac{n}{n+4}\dfrac{16\pi^4}{\big(B_nV(\Omega)\big)^{\frac{4}{n}}}k^{\frac{4}{n}}\\
 &+\left(\frac{n+2}{12n(n+4)}-\frac{1}{1152n^2(n+4)}
    \right)\frac{V(\Omega)}{I(\Omega)}\frac{n}{n+2}
 \dfrac{4\pi^2}{\big(B_nV(\Omega)\big)^{\frac{2}{n}}}k^{\frac{2}{n}}\\
 &+\left(\frac{1}{576n(n+4)}-\frac{1}{27648n^2(n+2)(n+4)}\right)
 \left(\frac{V(\Omega)}{I(\Omega)}\right)^2.
\endaligned
$$
This completes the proof of Theorem.
$$\eqno{\Box}$$

\vskip 3pt\noindent {\it Proof of  Corollary.}  Let $v_1$, $\cdots$,
$v_n$ be $n$ orthonormal eigenfunctions corresponding to the first $n$ eigenvalues $\mu_1$,
$\cdots$, $\mu_n$ of the Neumann eigenvalue problem of  Laplacian, that is,
\begin{equation*}
  {\begin{cases}
     \Delta v_i =- \mu_i  v_i,& \ \ {\rm in} \ \ \Omega ,\\[3mm]
     \dfrac{\partial v_i}{\partial \nu}=0 , & \ \ {\rm on}  \ \ \partial
     \Omega,\\[3mm]
     \int_\Omega v_iv_j=\delta_{ij},& \ \ {\rm i,j
     =1, \cdots, n}.
     \end{cases}}
\end{equation*}
 It then follows from the inequality (2.8) in \cite{[AB]} that
$$
\sum_{i=1}^n\frac{1}{\mu_i}\geq\frac{\int_\Omega |x|^2dx}{V(\Omega)}.\eqno{(2.29)}
$$
Combining (1.13) and (2.29), we have the inequality (1.14).
$$\eqno{\Box}$$

\end {document}